\documentclass[final]{siamltex}
\usepackage{amsmath}
\usepackage[square,comma,numbers,sort&compress]{natbib} %advanced bib package
\usepackage{eufrak}

\newtheorem{myremark}{Remark}[section]

\def\RR{\ensuremath{\bf R}}

\def\X{{\mathcal X}}
\def\Y{{\mathcal Y}}
\def\Z{{\mathcal Z}}
\def\H{{\mathcal H}}

\def\span{{\mbox{\rm span}}}

\def\lmax{\lambda_{\max}}
\def\lmin{\lambda_{\min}}
\def\dim{{\mbox{\rm dim}}}
\def\max{{\mbox{\rm max}}}
\def\min{{\mbox{\rm min}}}
\def\moo{\mathfrak{M}_{00}}
\def\moi{\mathfrak{M}_{01}}
\def\mio{\mathfrak{M}_{10}}
\def\mii{\mathfrak{M}_{11}}
\def\mo{\mathfrak{M}_{0}}

\title{Majorization for Changes in Angles Between Subspaces,
Ritz values, and graph Laplacian spectra
\thanks{Received by the editors August 31, 2005;
accepted for publication (in revised form) February 22, 2006;
published electronically \today.
Preliminary results of this paper were presented at the
12th ILAS Conference in Canada in June 2005
at the mini-symposium ``Canonical Angles Between Subspaces: Theory and Applications.''
}}

\author{Andrew V. Knyazev\thanks{Department of
Mathematical Sciences University of Colorado at Denver
and Health Sciences Center,
P.O. Box 173364, Campus Box 170, Denver, CO 80217-3364
(andrew.knyazev[at]cudenver.edu, http://math.cudenver.edu/\~{}aknyazev/).
This material is based upon work supported by the National Science
Foundation award DMS 0208773
and the Intelligence Technology Innovation Center through the joint
``Approaches to Combat Terrorism'' NSF Program Solicitation NSF 03-569.}
\and Merico E. Argentati\thanks{Department of
Mathematical Sciences University of Colorado at Denver
and Health Sciences Center,
P.O. Box 173364, Campus Box 170, Denver, CO 80217-3364 (rargenta[at]math.cudenver.edu).}}

\begin{document}
%\slugger{simax}{2006}{?}{?}{?--?}
\setcounter{page}{1}
\maketitle
\begin{abstract}
Many inequality relations between real vector quantities can be succinctly expressed as
``weak (sub)majorization'' relations using the symbol ${\prec}_{w}$.
We explain these ideas and apply them in several areas:
angles between subspaces, Ritz values, and graph Laplacian spectra,
which we show are all surprisingly related.

Let $\Theta(\X,\Y)$ be the vector of principal angles
in nondecreasing order
between subspaces $\X$ and $\Y$ of a finite dimensional space $\H$
with a scalar product. We consider the change in principal angles
between subspaces $\X$ and $\Z$, where we
let $\X$ be perturbed to give $\Y$. We measure the change using
weak majorization. We prove that
$|\cos^2\Theta(\X,\Z)-\cos^2\Theta(\Y,\Z)| {\prec}_{w} \sin\Theta(\X,\Y)$,
and give similar results for differences of cosines, i.e.
$|\cos\Theta(\X,\Z)-\cos\Theta(\Y,\Z)| {\prec}_{w} \sin\Theta(\X,\Y)$,
and of sines, and of sines squared, assuming $\dim \X = \dim \Y$.

We observe that $\cos^2\Theta(\X,\Z)$ can be interpreted
as a vector of Ritz values, where the Rayleigh-Ritz method is applied to the
orthogonal projector on $\Z$  using $\X$ as a trial subspace.
Thus, our result for the squares of cosines can be viewed
as a bound on the change in the Ritz values of an orthogonal projector.
We then extend it to prove a general result for Ritz values
for an arbitrary Hermitian operator $A$,
not necessarily a projector: let $\Lambda\left(P_{\X}A|_{\X}\right)$
be the vector of Ritz values in nonincreasing order for $A$
on a trial subspace $\X$,
which is perturbed to give another trial subspace $\Y$, then
$\left| \Lambda\left(P_{\X}A|_{\X}\right)
-  \Lambda\left(P_{\Y}A|_{\Y}\right) \right|
\prec_w (\lmax-\lmin)~\sin\Theta(\X,\Y)$,
where the constant is the difference between the largest and the
smallest eigenvalues of $A$.
This establishes our conjecture that the root two factor
in our earlier estimate may be eliminated.
Our present proof is based on a classical but rarely used technique of
extending a Hermitian operator in $\H$ to an orthogonal projector in
the ``double'' space $\H^2$.

An application of our Ritz values weak majorization result
for Laplacian graph spectra comparison is suggested, based
on the possibility to interpret eigenvalues of the edge Laplacian
of a given graph as Ritz values of the edge Laplacian of the
complete graph. We prove that
$ \sum_k |\lambda^1_k - \lambda^2_k| \leq n l,$
where $\lambda^1_k$ and $\lambda^2_k$ are all ordered
elements of the Laplacian spectra of two graphs
with the same $n$ vertices and with
$l$ equal to the number of differing edges.
\end{abstract}

\begin{AM}
15A42, %Inequalities involving eigenvalues and eigenvectors
15A60, %Norms of matrices, numerical range, applications of functional analysis to matrix theory
65F35, %Matrix norms, conditioning, scaling
05C50 %Graphs and matrices
\end{AM}
%\begin{DOI}
%?/?
%\end{DOI}

\pagestyle{myheadings}
\thispagestyle{plain}
\markboth{ANDREW V. KNYAZEV AND MERICO E. ARGENTATI}
{MAJORIZATION: ANGLES, RITZ VALUES, GRAPH SPECTRA}
%{MAJORIZATION FOR ANGLES, RITZ VALUES, AND GRAPH LAPLACIAN SPECTRA}

\section{Introduction}
Many inequality relations between real vector quantities can be
succinctly expressed as ``weak (sub)majorization'' relations using
the symbol ${\prec}_{w}$ that we now introduce.
For a real vector
$x = \left[ x_{1},\cdots,x_{n} \right]$
let $x^\downarrow$ be the
vector obtained by rearranging the entries of $x$ in an
algebraically non-increasing order.
Vector $y$ weakly majorizes vector $x$, i.e.
$x\prec_wy$,
if $\sum_{i=1}^k x_i^\downarrow \le \sum_{i=1}^k y_i^\downarrow, \ k=1,\ldots,n$.
The importance of weak majorization can
be seen from the classical statement that the following two conditions are equivalent:
$x{\prec}_{w}y$ and
$\sum _{i=1}^n \phi(x_i) \leq \sum_{i=1}^n \phi(y_i)$
for all nondecreasing convex functions $\phi$. Thus, a single weak
majorization result implies a great variety of inequalities. We
explain these ideas and apply them in several areas: angles
between subspaces, Ritz values, and graph Laplacian spectra, which
we show are all surprisingly related.

The concept of principal angles, also referred to as canonical angles,
between subspaces is one of the classical mathematical ideas
originated from \citet{MR1503705}
with many applications.
In functional analysis, the gap between subspaces, which is related to
the sine of the largest principal angle,
bounds the perturbation of a closed linear operator by measuring the change
in its graph, while the smallest nontrivial principal angle
between two subspaces determines if the sum of the subspaces is closed.
In numerical analysis, principal angles appear naturally to estimate
how close an approximate eigenspace is to the true eigenspace.
The chordal distance,
the Frobenius norm of the sine of the principal angles,
on the Grassmannian space of finite dimensional subspaces is used,
e.g., for subspace packing with applications in control theory.
In statistics,
the cosines of principal angles are called canonical correlations and
have applications in information retrieval and data visualization.

Let $\H$ be a real or complex $n < \infty$ dimensional vector space
equipped with an inner product
$(x,y)$  and a vector norm $\|x\|=(x,x)^{1/2}$.
The acute angle between two non--zero vectors
$x$ and $y$
is defined as
$$\nonumber
\theta({x},{y})=\arccos\frac{|(x,y)|}{\|x\|\|y\|} \in [0, \pi/2].
$$
For three nonzero vectors $x,y,z$, we have bounds on the
change in the angle
\begin{equation}\label{angle_only}
\left| \theta(x,z)-\theta(y,z)\right| \leq \theta(x,y),
\end{equation}
in the sine
\begin{equation}\label{angle_sines}
\left| \sin(\theta(x,z))-\sin(\theta(y,z))\right| \leq \sin(\theta(x,y)),
\end{equation}
in the cosine
\begin{equation}\label{angle_cosines}
\left| \cos(\theta(x,z))-\cos(\theta(y,z))\right| \leq \sin(\theta(x,y)),
\end{equation}
and a more subtle bound on the change in the sine or cosine squared
\begin{equation}\label{angle_cosine_squared}
\left| \cos^2(\theta(x,z))-\cos^2(\theta(y,z))\right|=
\left| \sin^2(\theta(x,z))-\sin^2(\theta(y,z))\right| \leq \sin(\theta(x,y)).
\end{equation}
Let us note that we can project the space $\H$ into the $\span \{x,y,z\}$ without
changing the angles, i.e. the inequalities above present
essentially the case of a $3D$ space.

Inequality (\ref{angle_only}) is proved in \citet[Th. 3.2, p. 514]{LiQiu01}.
We note that (\ref{angle_sines}) follows from (\ref{angle_only}),
since the sine function is increasing and subadditive, see \citet[p. 530]{LiQiu01}.

It is instructive to provide a simple proof of
the sine inequality (\ref{angle_sines}) using orthogonal projectors.
Let $P_\X$, $P_\Y$, and $P_\Z$ be, respectively, the
orthogonal projectors onto the subspaces spanned by the vectors $x$, $y$, and $z$,
and let $\|\cdot\|$ also denote the induced operator norm.
When we are dealing with $1D$ subspaces,
we have the following elementary formula
$\sin(\theta(x,y))=\|P_\X-P_\Y\|$
(indeed, $P_x-P_y$ has rank at most two, so it has at most
two non-zero singular values, but
$(P_x-P_y)^2x=(1-|(x,y)|^2)x$ and $(P_x-P_y)^2y=(1-|(x,y)|^2)y$
for unit vectors  $x$ and $y$, so
$1-|(x,y)|^2=\sin^2(\theta\{x,y\})$ is a double eigenvalue of
$(P_x-P_y)^2$).
Then the sine (\ref{angle_sines}) inequality is equivalent to
the triangle inequality
$\left|~\|P_\X - P_\Z\|-\|P_\Y-P_\Z\|~\right| \leq \|P_\X-P_\Y\|.$

In this paper, we replace
$1D$ subspaces spanned by the vectors $x$, $y$ and $z$,
with multi--dimensional subspaces  $\X$, $\Y$ and $\Z$,
and we use the concept of principal angles between subspaces.
Principal angles are very well studied in the literature, however, some
important gaps still remain. Here, we are
interested in generalizing inequalities
(\ref{angle_sines})-(\ref{angle_cosine_squared})
above to multi--dimensional subspaces
to include all principal angles, using weak majorization.

Let us denote by $\Theta(\X,\Y)$ the vector of principal angles
in nondecreasing order between subspaces $\X$ and $\Y$.
Let $\dim \X = \dim \Y$, and let another subspace $\Z$ be given.
We prove that
$|\cos^2\Theta(\X,\Z)-\cos^2\Theta(\Y,\Z)| {\prec}_{w} \sin\Theta(\X,\Y)$,
and give similar results for differences of cosines, i.e.
$|\cos\Theta(\X,\Z)-\cos\Theta(\Y,\Z)| {\prec}_{w} \sin\Theta(\X,\Y)$,
and of sines, and of sines squared.
This is the first main result of the present paper,
see Section \ref{section_angles_proximity}.
The proof of weak majorization for sines is a direct
generalization of the 1D proof above. Our proofs of
weak majorization for cosines and sines or cosines squared
do not have such simple 1D analogs.

Pioneering results using angles between
subspaces in the framework of unitarily invariant norms and
symmetric gauge functions, equivalent to majorization,
appear in \citet{davis1}, which introduces many of the tools that we use here.
The main goal of \citet{davis1} is however different ---
analyzing the perturbations of eigenvalues and eigenspaces,
while in the present paper
we are concerned with sensitivity of angles and
Ritz values with respect to changes in subspaces.

Our second main result, see Section \ref{section_ritz},
bounds the change in the Ritz values with the change of the trial
subspace. We attack the problem by discovering a simple, but deep, connection
between the principal angles and the Rayleigh--Ritz method.

We first give a brief definition of Ritz values.
Let $A:\H \to\H$ be a Hermitian operator and
let ${\X}$ be a (so-called ``trial'') subspace of $\H$.
We define an operator
$P_{\X}A|_{\X}$ on $\X$, where $P_{\X}$ is the orthogonal projector
onto $\X$ and $P_{\X}A|_{\X}$ denotes the restriction of
operator $P_{\X}A$ to its invariant subspace $\X$,
as discussed, e.g., in  \citet[Sect. 11.4, pp. 234--239]{p}.
The eigenvalues $\Lambda(P_{\X}A|_{\X})$ are called Ritz values
of the operator $A$ with respect to the trial subspace $\X$.

We observe that the cosines squared $\cos^2\Theta(\X,\Z)$
of principal angles
between subspaces $\X$ and $\Z$ can be interpreted
as a vector of Ritz values, where the Rayleigh-Ritz method is applied to the
orthogonal projector $P_\Z$ onto $\Z$  using $\X$ as a trial subspace.
Let us illustrate this connection
for one--dimensional $\X=\span\{x\}$ and $\Z=\span\{z\}$,
where it becomes trivial:
$$
\cos^2(\theta(x,z)) = \frac{(x,P_\Z x)}{(x,x)}.
$$
The ratio on the right is the Rayleigh quotient for $P_\Z$ ---
the one dimensional analog of a Ritz value.
In this notation, estimate (\ref{angle_cosine_squared}) turns into
\begin{equation}\label{angle_cosine_squared_RR}
\left| \frac{(x,P_z x)}{(x,x)} - \frac{(y,P_z y)}{(y,y)} \right| \leq \sin(\theta(x,y)),
\end{equation}
which clearly now is a particular case of a
general estimate for the Rayleigh quotient, cf. \citet{ka03},
\begin{equation}\label{angle_RR}
\left| \frac{(x,A x)}{(x,x)} - \frac{(y,A y)}{(y,y)} \right| \leq
\left(\lambda_\max - \lambda_\min \right) \sin(\theta(x,y)),
\end{equation}
where $A$ is a Hermitian operator and $\lambda_\max - \lambda_\min$ is the spread
of its spectrum.

We show that the multi--dimensional analog of (\ref{angle_cosine_squared_RR})
can be interpreted as a bound on the
change in the Ritz values with the change of the trial
subspace, in the particular case where the Rayleigh-Ritz method is applied to an
orthogonal projector.
We then extend it to prove a general result for Ritz values
for an arbitrary Hermitian operator $A$,
not necessarily a projector: let $\Lambda\left(P_{\X}A|_{\X}\right)$
be the vector of Ritz values in nonincreasing order for the operator $A$
on a trial subspace $\X$,
which is perturbed to give another trial subspace $\Y$, then
$\left| \Lambda\left(P_{\X}A|_{\X}\right)
-  \Lambda\left(P_{\Y}A|_{\Y}\right) \right|
\prec_w (\lmax-\lmin)~\sin\Theta(\X,\Y)$,
which is a multi--dimensional analog of (\ref{angle_RR}).
Our present proof is based on a classical but rarely used idea of
extending a Hermitian operator in $\H$ to an orthogonal projector in
the ``double'' space $\H^2$, preserving its Ritz values.

An application of our Ritz values weak majorization result
for Laplacian graph spectra comparison is suggested
in Section \ref{s.g}, based
on the possibility to interpret eigenvalues of the edge Laplacian
of a given graph as Ritz values of the edge Laplacian of the
complete graph. We prove that
$ \sum_k |\lambda^1_k - \lambda^2_k| \leq n l,$
where $\lambda^1_k$ and $\lambda^2_k$ are all ordered
elements of the Laplacian spectra of two graphs
with the same $n$ vertices and with
$l$ equal to the number of differing edges.

The rest of the paper is organized as follows.
In Section \ref{section_definitions}, we provide
some background, definitions and several statements concerning
weak majorization,
principal angles between subspaces,
and extensions of Hermitian operators to projectors.
In Section \ref{section_angles_proximity},
we prove
in Theorems \ref{thm.sine} and \ref{thm.cosine_squared}
that the absolute value of the
change in (the squares of) the sines and cosines
is weakly majorized by the sines of the angles between
the original and perturbed subspaces.
In Section \ref{section_ritz}, we prove
in Theorem \ref{thm_ritz}
that a change in the Ritz values
in the Rayleigh-Ritz method
with respect to the change in the trial subspaces is
weakly majorized by the sines of the principal angles between
the original and perturbed trial subspaces times a constant.
In Section \ref{s.g} we apply our Ritz values weak majorization result
to Laplacian graphs spectra comparison.

This paper is related to several different subjects: majorization,
principal angles, Rayleigh-Ritz method, and Laplacian graph
spectra. In most cases, whenever possible, we cite books rather
than the original works in order to keep our already quite long
list of references within a reasonable size.

\section{Definitions and Preliminaries}\label{section_definitions}
In this section we introduce some definitions,
basic concepts and mostly known results for later use.

\subsection{Weak Majorization}\label{sect.majorization}
Majorization is a well known,
e.g.,\ in \citet[pp. 45--49]{Hardy:1934:I} or \citet[pp. 9--14]{mo},
important mathematical concept with numerous applications.

For a real vector
$x = \left[ x_{1},\cdots,x_{n} \right]$
let $x^\downarrow$ be the
vector obtained by rearranging the entries of $x$ in an
algebraically non-increasing order,
$x_{1}^\downarrow \geq \cdots \geq x_{n}^\downarrow.$
We denote $[|x_1|,\cdots,|x_n|]$ by $|x|$.
We say that vector $y$ weakly majorizes vector $x$ and we use the notation
$[x_1,\cdots,x_n] \prec_w [y_1,\cdots,y_n]$
or $x\prec_wy$~~
if ~~ $\sum_{i=1}^k x_i^\downarrow \le \sum_{i=1}^k y_i^\downarrow,\quad \ k=1,\ldots,n$.
If in addition the sums above for $k=n$ are equal,
 $y$ (strongly) majorizes vector $x$, but we do not use
this type of majorization in the present paper.
Two vectors of different lengths may be compared by simply appending zeroes to
increase the size of the smaller vector to make the vectors the same length.

Weak majorization is a powerful tool for estimates involving
eigenvalues and singular values and is covered, e.g.,\ in
\citet{gohKr}, \citet{mo}, \citet{bhatia_book} and \citet{hj2}, which
we follow here and refer the reader to
for references to the original works and all necessary proofs.
In the present paper, we use several well known statements
that we formulate for operators $\H\to\H$ and
overview briefly below.

Let $S(A)$ denote the vector of all singular values of $A:\H\to\H$ in
nonincreasing order, i.e. $S(A)=S^\downarrow(A)$,
while individual singular values of $A$ enumerated in
nonincreasing order are denoted by $s_i(A)$.
For Hermitian $A$ let $\Lambda(A)$ denote the vector of all eigenvalues of $A$
in nonincreasing order, i.e. $\Lambda(A)=\Lambda^\downarrow(A)$,
while individual eigenvalues of $A$ enumerated in
nonincreasing order are denoted by $\lambda_i(A)$.

The starting point for
weak majorization results we use in this paper is,
e.g., \citet[Th. 9.G.1, p. 241]{mo},
\begin{theorem} \label{thm.e_sum_majorization}
 $\Lambda(A+B) \prec_w \Lambda(A) + \Lambda(B)$ for Hermitian $A$ and $B$.
\end{theorem}\\
This follows easily from Ky Fan's trace maximum principle
\citet[Th. 20.A.2, p. 511]{mo}
and the fact that the maximum of a sum is bounded from above by the
sum of the maxima.
For general $A:\H\to\H$ and $B:\H\to\H$, it follows from Theorem \ref{thm.e_sum_majorization},
since the top half of the spectrum of the Hermitian 2-by-2 block operator
$\left[\begin{array}{cc}
0 & A \\
A^\ast & 0 \end{array} \right]$
is nothing but $S(A)$,
that, e.g., \citet[Cor. 3.4.3, p. 196]{hj2},
\begin{corollary} \label{thm.s_sum_majorization}
$S(A\pm B) \prec_w S(A) + S(B)$.
\end{corollary}\\
A more delicate and stronger result is the following Lidskii theorem e.g.,
\citet[Th. III.4.1, p. 69 ]{bhatia_book},
which can be proved using the Wielandt maximum principle,
e.g., \citet[Th. III.3.5, p. 67]{bhatia_book},
\begin{theorem} \label{thm.e_sum_majorization_s}
For Hermitian $A$ and $B$ and any set
of indices $1\leq i_1<\ldots<i_k\leq n=\dim \H$, we have
$
\sum_{j=1}^k \lambda_{i_j}(A+B)
\leq
\sum_{j=1}^k \lambda_{i_j}(A) +
\sum_{j=1}^k \lambda_{j}(B), \, k=1,\ldots,n.
$
\end{theorem}
By choosing an appropriate set of indices, Theorem
\ref{thm.e_sum_majorization_s} for Hermitian $A$ and $B$ immediately
gives $\Lambda(A) - \Lambda(B)  \prec_w \Lambda(A-B)$, which for
singular values of arbitrary $A:\H\to\H$ and $B:\H\to\H$ is
equivalent (see \citet[Sect. IV.3, pp. 98--101]{bhatia_book}), to,
e.g. \citet[Th. 3.4.5, p. 198]{hj2} or \citet[Th. IV.3.4, p. 100]{bhatia_book},
\begin{corollary}\label{cor.s_dif_majorization}
$\left| S(A) - S(B) \right|  \prec_w S(A-B)$.
\end{corollary}\\
Applying Corollary \ref{cor.s_dif_majorization}
to properly shifted Hermitian operators, we get
\begin{corollary}\label{cor.e_dif_majorization}
$\left| \Lambda(A) - \Lambda(B) \right|  \prec_w S(A-B)$ for Hermitian $A$ and $B$.
\end{corollary}\\
We finally need the so called ``pinching'' inequality,
e.g. \citet[Th. II.5.1, p. 52]{gohKr} or \citet[Prob. II.5.5, p. 50]{bhatia_book},
\begin{theorem} \label{thm.pinching}
If $P$ is an orthogonal projector then
$S(PAP\pm(I-P)A(I-P)) \prec_w S(A).$
\end{theorem}
\begin{proof}
Indeed, $A=PAP+(I-P)A(I-P)+PA(I-P)+(I-P)AP$ so
let $B=PAP+(I-P)A(I-P)-PA(I-P)-(I-P)AP$ then
$(2P-I)A(2P-I)=B$, where $2P-I$ is unitary Hermitian, so
$A^\ast A$ and $B^\ast B$ are similar and $S(A)=S(B)$. Evidently,
$PAP+(I-P)A(I-P) = (A+B)/2$  so the pinching result with the plus
follows from Corollary \ref{thm.s_sum_majorization}.
The pinching result with the minus
is equivalent to the pinching result with the plus
since the sign does not change the singular values on the
left-hand side:
$S(PAP\pm(I-P)A(I-P)) = S(PAP) \cup S((I-P)A(I-P)),$
since the ranges of $PAP$ and $(I-P)A(I-P)$ are disjoint.
\end{proof}

\subsection{Principal Angles Between Subspaces}
Let $P_{\X}$ and $P_{\Y}$ be orthogonal projectors onto the %(nontrivial proper)
subspaces $\X$ and $\Y$, respectively, of the space $\H$.
We define the set of cosines of principal angles between
 subspaces $\X$ and $\Y$ by
\begin{equation}\label{eq.def_angles}
\cos \Theta(\X,\Y)= [s_1(P_{\X}P_{\Y}),\ldots,s_m(P_{\X}P_{\Y})],\,
m = \min\left\{\dim \X; \dim \Y\right\}.
\end{equation}
Our definition (\ref{eq.def_angles}) is evidently symmetric:
$\Theta(\X,\Y) = \Theta(\Y,\X).$
By definition, the cosines are arranged in nonincreasing order,
i.e. $\cos(\Theta(\X,\Y)) = (\cos(\Theta(\X,\Y)))^\downarrow,$
while the angles $\theta_i(\X,\Y) \in [0,{\pi}/{2}], i=1,\ldots,m$
and their sines are in nondecreasing order.

The concept of principal angles is closely
connected to cosine--sine (CS) decompositions of
unitary operators; and we refer the reader to the books
\citet{ss,s2,bhatia_book}
for the history and
references to the original publications on the principal angles
and the CS decomposition.
We need several simple but important statements about the angles
provided below.
In the particular case $\dim \X = \dim \Y$,
the standard CS decomposition can be used and
the statements are easy to derive.
For the general case $\dim \X \neq \dim \Y$ that is necessary for us here,
they can be obtained using the general (rectangular) form
of the CS decomposition
described, e.g., in \citet{MR95f:15001,MR615522}.
For completeness we provide the proofs here
using ideas from \citet{MR0251519,davis1},
preparing our work to be more easily extended
to infinite dimensional Hilbert spaces.

We formulate without proofs in \citet[Th. 3.4, p. 2017 and Th. 3.5, p. 2018]{ka02}
statements, equivalent to the theorem below, see also \citet[Ex. 1.2.6]{MR94d:65002}.
\begin{theorem}\label{thm.angles.relations}
When one of the two subspaces is replaced with its orthogonal
complement, the corresponding pairs of angles sum up to $\pi/2$,
specifically:
\begin{equation}\label{eq.perp1}
\left[\frac{\pi}{2},\ldots,\frac{\pi}{2}, (\Theta(\X,\Y))^\downarrow\right]=
\left[\frac{\pi}{2} - \Theta(\X,\Y^\perp), 0, \ldots, 0\right],
\end{equation}
where there are $\max\{\dim \X - \dim \Y;0\}$  values $\pi/2$ on the left,
and possibly extra zeros on the right to match the sizes.

The angles between subspaces and between
their orthogonal complements are essentially the same,
\begin{equation}\label{eq.perp}
[(\Theta(\X,\Y))^\downarrow, 0, \ldots, 0]= [(\Theta(\X^\perp,\Y^\perp))^\downarrow,0, \ldots, 0],
\end{equation}
where extra $0$s at the end may need to be added on either side to match the sizes.
\end{theorem}
\begin{proof}
Let
$
\moo=\X \cap \Y, \enspace \moi=\X \cap \Y^\perp, \enspace \mio=\X^\perp \cap \Y, \enspace
\mii=\X^\perp \cap \Y^\perp,
$
as suggested in \citet[p. 381]{MR0251519}.

Each of the
subspaces is invariant with respect to
orthoprojectors $P_{\X}$ and $P_{\Y}$ and their products,
and so each of the subspaces contributes independently
to the set of singular values of $P_{\X}P_{\Y}$ in
(\ref{eq.def_angles}). Specifically, there are $\dim \moo$
ones,  $\dim \mo$ singular values in the interval $(0,1)$
equal to $\cos \Theta(\mo,\Y)$,
where $\mo= \X \cap (\moo \oplus \moi)^\perp$,
and all other singular values are zeros; thus,
\begin{equation}\label{eq.angle.1}
\Theta(\X,\Y)^\downarrow=\left[\frac{\pi}{2},\ldots,\frac{\pi}{2},
(\Theta(\mo,\Y))^\downarrow, 0, \ldots, 0\right],
\end{equation}
where there are $\min\left\{\dim(\moi);\dim(\mio) \right\}$ values $\pi/2$ %on the left,
and $\dim(\moo)$ zeros. %on the right
%and we use the fact that $\dim \moi - \dim \mio = \dim X - \dim Y$.

The subspace $\mo$ does not change if we substitute $\Y^\perp$ for $\Y$ in (\ref{eq.angle.1}),
so we have
$$
\Theta(\X,\Y^\perp)^\downarrow=\left[\frac{\pi}{2},\ldots,\frac{\pi}{2},
(\Theta(\mo,\Y^\perp))^\downarrow, 0, \ldots, 0\right],
$$
where there are
$\min\left\{\dim(\moo);\dim(\mii) \right\}$ values $\pi/2$ %on the left,
and $\dim(\moi)$ zeros. %on the right.
Since $\lambda$ is an eigenvalue of $(P_\X P_\Y) |_{\mo}$ if and only if $1-\lambda$ is
an eigenvalue of $(P_\X P_{\Y^\perp}) |_{\mo}$, we have
$\frac{\pi}{2} - \Theta(\mo,\Y^\perp)=(\Theta(\mo,\Y))^\downarrow$,
and the latter equality turns into
\begin{equation}\label{eq.angle.2}
\frac{\pi}{2} - \Theta(\X,\Y^\perp)=\left[\frac{\pi}{2},\ldots,\frac{\pi}{2},
(\Theta(\mo,\Y))^\downarrow, 0, \ldots, 0\right],
\end{equation}
where there are $\dim(\moi)$  values $\pi/2$, %on the left,
and $\min\left\{\dim(\moo);\dim(\mii) \right\}$ zeros. %on the right.
To obtain (\ref{eq.perp1}), we make %the right hand sides, and thus the left hand sides, of
(\ref{eq.angle.1}) and (\ref{eq.angle.2}) equal by
adding $\max\left\{\dim \moi - \dim \mio;0 \right\}$ values $\pi/2$ to
(\ref{eq.angle.1}) and $\max\left\{\dim \moo - \dim \mii;0 \right\}$ zeros to
(\ref{eq.angle.2}), and
noting that, since
$\dim(\X \cap \Y^\perp)=\dim\X + \dim\Y^\perp - \dim(\X+\Y^\perp)$,
$\X^\perp \cap \Y=(\X+\Y^\perp)^\perp$, and
$\dim\Y^\perp - \dim((\X+\Y^\perp)^\perp) = \dim(\X+\Y^\perp)-\dim\Y$,
we have
$\dim \moi - \dim \mio = \dim(\X \cap \Y^\perp) - \dim(\X^\perp \cap \Y) =
\dim\X + \dim\Y^\perp - \dim(\X+\Y^\perp) - \dim((\X+\Y^\perp)^\perp) =
\dim\X - \dim\Y$.

The proof above shows that there are $\dim\moo= \dim(\X \cap \Y)$
zeros on the right in (\ref{eq.perp1}).
To prove (\ref{eq.perp}), we substitute in (\ref{eq.perp1})
$\X^\perp$ for $\X$ to get
$
\left[\frac{\pi}{2},\ldots,\frac{\pi}{2}, (\Theta(\X^\perp,\Y))^\downarrow\right]=
\left[\frac{\pi}{2} - \Theta(\X^\perp,\Y^\perp), 0, \ldots, 0\right]
$
with $\dim(\X^\perp \cap \Y)$ zeros on the right
on the one hand and exchange $\Theta(\X,\Y)=\Theta(\Y,\X)$ in (\ref{eq.perp1})
and then substitute $\X^\perp$ for $\X$ to obtain
$
\left[\frac{\pi}{2},\ldots,\frac{\pi}{2}, (\Theta(\Y,\X^\perp))^\downarrow\right]=
\left[\frac{\pi}{2} - \Theta(\Y,\X), 0, \ldots, 0\right]
$
with $\dim(\Y \cap \X^\perp)$ zeros on the right
on the other hand. We have equal numbers of zeros on the right
in both equalities and $\Theta(\X^\perp,\Y)=\Theta(\Y,\X^\perp)$
by the symmetry of our definition (\ref{eq.def_angles}),
so subtracting both equalities from $\pi/2$ leads to (\ref{eq.perp}).
\end{proof}

We also use the following trivial, but crucial, statement.
\begin{lemma} \label{lem_projector_angles_pre}
$
\Lambda\left( (P_{\X} P_{\Y})|_{\X} \right)=[ \cos^2 \Theta(\X,\Y),0,\ldots,0],
$
with $\max\{\dim\X - \dim\Y;0\}$ extra $0$s.
\end{lemma}
\begin{proof}
The operator
$(P_{\X} P_{\Y})|_{\X}=\left((P_{\X} P_{\Y}) (P_{\X} P_{\Y})^{\star}\right)|_{\X}$
is Hermitian nonnegative definite,
and its spectrum can be represented using
the definition of angles (\ref{eq.def_angles}).
The number of extra $0$s is exactly the difference between the
number $\dim\X$ of  Ritz values and the
number $\min\{\dim\X;\dim\Y\}$ of principal angles .
\end{proof}

Finally, we need the following characterization of
singular values of the difference of projectors, which for $\dim\X = \dim\Y$ appears, e.g.,\ in
\citet[Th. 5.5.5, p. 43]{s2}.
\begin{theorem}\label{th.SFmG}
$$
[S(P_{\X}-P_{\Y}), 0, \ldots, 0]
= [1,\ldots,1,\left(\sin \Theta(\X,\Y),\sin \Theta(\X,\Y)\right)^\downarrow, 0, \ldots, 0],
%\end{equation}
$$
where there are $|\dim \X - \dim \Y |$ extra $1$s upfront, the set
$\sin \Theta(\X,\Y)$ is repeated twice and ordered,
and extra $0$s at the end may need to be added
on either side to match the sizes.
\end{theorem}
\begin{proof}
The projectors $P_\X$ and $P_\Y$ are idempotent, which implies on the one hand
$$
(P_\X-P_\Y)^2 = P_\X (I- P_{\Y}) + P_\Y (I-P_{\X})=
P_\X P_{\Y^\perp} + P_\Y P_{\X^\perp},
$$
so the subspace $\X$ is invariant under $(P_\X-P_\Y)^2.$
On the other hand,
$$
(P_\X-P_\Y)^2 = (I- P_{\X}) P_\Y  + (I-P_{\Y}) P_\X= P_{\X^\perp} P_\Y + P_{\Y^\perp} P_\X,
$$
so the
subspace $\X^\perp$ is also invariant under $(P_\X-P_\Y)^2$.
The projectors $P_\X$ and $P_\Y$ are orthogonal, thus
the operator $(P_\X-P_\Y)^2$ is Hermitian, and
its spectrum can be represented as a union (counting the multiplicities)
of the spectra of its restrictions to the
complementary invariant subspaces $\X$ and $\X^\perp$:
$$
\Lambda \left((P_{\X}-P_{\Y})^2\right)=\left[
\Lambda((P_\X P_{\Y^\perp})|_\X),\Lambda((P_{\X^\perp} P_\Y))|_{\X^\perp}\right]^\downarrow .
$$
Using Lemma \ref{lem_projector_angles_pre} and statement (\ref{eq.perp1}) of
Theorem \ref{thm.angles.relations},
\begin{eqnarray}\nonumber
[\Lambda((P_\X P_{\Y^\perp})|_\X), 0, \ldots, 0]
&=&[\cos^2 \Theta(\X,\Y^\perp), 0, \ldots, 0]\\\nonumber
&=&[1,\ldots,1,(\sin^2 \Theta(\X,\Y))^\downarrow,0, \ldots, 0],
\end{eqnarray}
where there are $\max\{\dim \X - \dim \Y;0\}$  leading $1$s %on the right
and possibly extra zeros to match the sizes, and
\begin{eqnarray}\nonumber
[\Lambda((P_{\X^\perp} P_\Y))|_{\X^\perp}), 0, \ldots, 0]
&=&[\cos^2 \Theta(\X^\perp,\Y), 0, \ldots, 0]\\\nonumber
&=&[1,\ldots,1,(\sin^2 \Theta(\X,\Y))^\downarrow,0, \ldots, 0],
\end{eqnarray}
where there are $\max\{\dim \Y - \dim \X;0\}$  leading $1$s %on the right
and possibly extra zeros to match the sizes.
Combining these two relations
and taking the square root completes the proof.
\end{proof}

\subsection{Extending Operators to Isometries and Projectors}\label{section_extend_proj}
In this subsection we present a simple and known technique,
e.g.,\ \citet{h50} and \citet[p. 461]{rsn}
for extending a Hermitian operator to a projector.
We give an alternative proof based on extending an arbitrary normalized operator $B$ to
an isometry $\hat{B}$ (in matrix terms, a matrix with orthonormal columns).
\citet[Prob. X.1.26, p. 455]{MR0354718} and \citet[Prob. I.3.6, p. 11]{bhatia_book} extend
$B$ to a block $2$-by-$2$ unitary operator.
Our technique is similar and results in a
$2$-by-$1$ isometry operator $\hat{B}$ that
coincides with the first column of the
$2$-by-$2$ unitary extension.
\begin{lemma} \label{isometric_extension}
Given an operator $B:\H\to\H$ with singular values less
than or equal to one, there exists a block $2$-by-$1$ isometry operator
$\hat{B}:\H\to\H^2$,
such that the upper block of $\hat{B}$ coincides with $B$.
\end{lemma}
\begin{proof}
$B^{\ast}B$ is Hermitian nonnegative definite,
and all its eigenvalues are bounded by one,
since all singular values of $B$ are bounded by one.
Therefore, $I-B^{\ast}B$ is Hermitian and nonnegative definite,
and thus possesses a Hermitian nonnegative square root.
Let
$$
\hat{B}=\left[\begin{array}{c}
B\\
\sqrt{I-B^{\ast}B}
\end{array} \right].
$$
By direct calculation,
$
\hat{B}^\ast\hat{B}=B^\ast B+\sqrt{I-B^{\ast}B}\sqrt{I-B^{\ast}B}=I,
$
i.e. $\hat{B}$ is an isometry.
\end{proof}

Now we use Lemma \ref{isometric_extension} to extend,
in a similar sense, a shifted and normalized Hermitian
operator to an orthogonal projector.
\begin{theorem} (\citet{h50} and \citet[p. 461]{rsn})
\label{thm.ext_to_projector}
Given a Hermitian operator $A:\H\to\H$ with eigenvalues enclosed
in the segment $\left[ 0,1\right]$, there exists a  block $2$-by-$2$
orthogonal projector $\hat{A}:\H^2\to\H^2$,
such that its upper left block is equal to $A$.
\end{theorem}
\begin{proof}
There exists $\sqrt{A}$,
which is also Hermitian and has its eigenvalues enclosed in $\left[ 0,1\right]$. Applying
Lemma \ref{isometric_extension} to $B=\sqrt{A}$, we construct the isometry $\hat{B}$
and set
$$
\hat{A}=\hat{B}\hat{B}^\ast=\left[\begin{array}{c}
\sqrt{A}\\
\sqrt{I-A}
\end{array} \right]\left[\begin{array}{cc}
\sqrt{A} & \sqrt{I-A}
\end{array} \right]=\left[\begin{array}{cc}
A & \sqrt{A(I-A)} \\
\sqrt{A(I-A)} & I-A \end{array} \right].
$$
We see that indeed the upper left block is equal to $A$.
We can use the fact that $\hat{B}$ is an isometry to show that
 $\hat{A}$ is an orthogonal projector, or that can be checked
directly by calculating $\hat{A}^2=\hat{A}$ and noticing that
$\hat{A}$ is Hermitian by construction.
\end{proof}

Introducing  $S=\sqrt{A}$ and $C=\sqrt{I-A}$, we obtain
$$
\hat{A}=\left[\begin{array}{cc}
S^2 & SC \\
SC & C^2 \end{array} \right],
$$
which is a well known, e.g., \citet{MR0251519,MR0098980},
block form of an orthogonal projector that can alternatively be
derived using the CS decomposition of unitary operators,
e.g., \citet{ss,s2,bhatia_book}.

The importance of Theorem \ref{thm.ext_to_projector}
can be better seen if we reformulate it as
\begin{theorem}
\label{thm.A_to_angles}
(cf. \citet[Ex. X.27, p. 455]{MR0354718})
Given a Hermitian operator $A:\H\to\H$ with eigenvalues enclosed
in a segment $\left[ 0,1\right]$, there exist subspaces
$\X$ and $\Y$ in $\H^2$ such that $A$
is unitarily equivalent to
$\left.(P_\X P_\Y)\right|_\X,$
where $P_\X$ and $P_\Y$ are the corresponding orthogonal
projectors in $\H^2$ and $|_\X$ denotes a restriction
to the invariant subspace $\X$.
\end{theorem}
\begin{proof}
We  use Theorem \ref{thm.ext_to_projector}
and take $P_\Y=\hat{A}$ and $P_\X = \left[\begin{array}{cc}
I & 0 \\
0 & 0 \end{array} \right]$.
\end{proof}

Similar to Theorem \ref{thm.A_to_angles},
Lemma \ref{isometric_extension} implicitly states that
an arbitrary normalized operator $B$ is unitary equivalent
to a product of the partial isometry $\hat{B}$  in $\H^2$
and the orthogonal projector in $\H^2$ that selects
the upper block in $\hat{B}$ (called $P_\X$ in the proof of
Theorem \ref{thm.A_to_angles}).
It is instructive to compare this product to
the classical polar decomposition of $B$ that is a product of a  partial isometry
and a Hermitian nonnegative operator in $\H$. In $\H^2$,
we can choose the second factor to be an orthogonal projector!
This statement together with Theorem \ref{thm.A_to_angles}
can provide interesting canonical decompositions in $\H^2$ that
apparently are not exploited at present, but in our opinion
deserve attention.

We take advantage of Theorem \ref{thm.A_to_angles} in the present paper.
Using Lemma \ref{lem_projector_angles_pre} with (\ref{eq.def_angles}),
Theorem \ref{thm.A_to_angles} implies that {\em the spectrum of an
arbitrary Hermitian operator after a proper shift and scaling is nothing
but a set of cosines squared of principal angles between some pair
of subspaces}. This surprising idea appears to be very powerful. It allows us,
in Section \ref{section_ritz}, to obtain a novel
result on sensitivity of Ritz values
with respect to the trial subspace by reducing the
investigation of the Rayleigh-Ritz method to the analysis
of the principal angles between subspaces that we provide
in the next section.

\section{Majorization for Angles}\label{section_angles_proximity}
In this section we prove the main results of the present paper
involving sines and cosines and their squares of principal angles,
but we start with
a known statement that involves the principal angles themselves:
\begin{theorem}(\citet[Th. 3.2, p. 514]{LiQiu01})\label{thm.Qiu}
Let ${\X}$, ${\Y}$ and $\Z$ be subspaces of
the same dimension. Then
\begin{equation}\label{eq.theta}
|\Theta(\X,\Z)-\Theta(\Y,\Z)| \prec_w \Theta(\X,\Y).
\end{equation}
\end{theorem}

Theorem \ref{thm.Qiu} deals with the principal angles themselves, and the obvious
question is: are there similar results for a function of these angles,
in particular for sines and cosines and their squares?
For one dimensional subspaces, estimate (\ref{eq.theta})
turns into (\ref{angle_only}) that, as discussed in the Introduction,
implies the estimate (\ref{angle_sines}) for the sine.
According to an anonymous referee,
it appears to be known to some specialists
that the same inference can be made for tuples of angles, but
there is no good reference for this at present.
Below we give easy direct proofs in a unified way for the
sines and cosines and their squares.

We first prove the estimates for sine and cosine, which are
straightforward generalizations of the 1D sine (\ref{angle_sines})
and cosine (\ref{angle_cosines}) inequalities from the Introduction.
\begin{theorem} \label{thm.sine}
Let $\dim \X = \dim \Y$ then
\begin{equation}\label{eq.sine_pert}
|\sin\Theta(\X,\Z)-\sin\Theta(\Y,\Z)| \prec_w \sin\Theta(\X,\Y),
\end{equation}
\begin{equation}\label{eq.cosine_pert}
|\cos\Theta(\X,\Z)-\cos\Theta(\Y,\Z)| \prec_w \sin\Theta(\X,\Y).
\end{equation}
\end{theorem}
\begin{proof}
Let $P_{\X}$, $P_{\Y}$ and $P_{\Z}$ be the corresponding orthogonal projectors onto
the subspaces $\X$, $\Y$ and $\Z$, respectively.
We prove the sine estimate (\ref{eq.sine_pert}), using the idea of \citet{LiQiu_05}.
Starting with
$
(P_\X - P_\Z ) - ( P_\Y - P_\Z) = P_\X - P_\Y,
$
as in the proof of the 1D sine estimate (\ref{angle_sines}),
we use Corollary \ref{cor.s_dif_majorization} to obtain
$$
| S(P_\X - P_\Z ) - S( P_\Y-P_\Z)|
\prec_w
S(P_\X - P_\Y).
$$
The singular values of the difference of two orthoprojectors are
described by Theorem \ref{th.SFmG}. Since $\dim \X = \dim \Y$ we have
the same number of extra $1$s upfront in  $S(P_\X - P_\Z )$ and in
$S(P_\Y - P_\Z )$ so that the extra $1$'s are cancelled and the
set of nonzero entries of $| S(P_\X - P_\Z ) - S( P_\Y-P_\Z)|$
consists of nonzero entries of $|\sin \Theta(\X,\Z) - \sin
\Theta(\Y,\Z)|$ repeated twice. The nonzero entries of $S(P_\X -
P_\Y)$ are by Theorem \ref{th.SFmG} the nonzero entries of $\sin
\Theta(\X,\Y)$ also repeated twice, thus we come to
(\ref{eq.sine_pert}).

The cosine estimate (\ref{eq.cosine_pert})
follows directly from the sine estimate
(\ref{eq.sine_pert}) with $\Z^\perp$ instead of $\Z$
because of (\ref{eq.perp1}) utilizing the assumption $\dim \X = \dim \Y$ .
\end{proof}

In our earlier paper, \citet[Lem. 5.1, p. 2023 and Lem. 5.2, p. 2025]{ka02}
we obtained a particular case of Theorem \ref{thm.sine},
only for the largest change in the sine and the cosine,
but with improved constants.
We are not presently able, however, to modify the proofs of \citet{ka02}
using weak majorization, in order to improve the
estimates of Theorem \ref{thm.sine}
by introducing these same constants.

Our last, but not least, result in this series is the weak majorization
bound for the sines or cosines squared, which provides the
foundation for the rest of the paper.
\begin{theorem} \label{thm.cosine_squared}
Let $\dim \X = \dim \Y$, then
%\begin{equation}\label{eq1.main_inequality}
$$
|\cos^2\Theta(\X,\Z)-\cos^2\Theta(\Y,\Z)| =
|\sin^2\Theta(\X,\Z)-\sin^2\Theta(\Y,\Z)|  \prec_w \sin\Theta(\X,\Y).
$$
%\end{equation}
%\begin{equation}\label{eq2.main_inequality}
%|\cos^2\Theta(\X,\Z)-\cos^2\Theta(\Y,\Z))| \prec_w \sin\Theta(\X,\Y).
%\end{equation}
\end{theorem}
\begin{proof}
The equality is evident.
To prove the majorization result for the sines squared,
we start with the useful pinching identity
$$
(P_\X -  P_\Z )^2 - ( P_\Y - P_\Z)^2 =
P_{\Z^\perp}(P_\X - P_\Y)P_{\Z^\perp} - P_\Z(P_\X - P_\Y)P_\Z.
$$
Applying Corollary \ref{cor.s_dif_majorization} we obtain
$$
\left| S\left((P_\X -  P_\Z )^2\right) - S\left(( P_\Y - P_\Z)^2\right) \right|
\prec_w
S\left(P_{\Z^\perp}(P_\X - P_\Y)P_{\Z^\perp} - P_\Z(P_\X - P_\Y)P_\Z \right).
$$
For the left-hand side we use Theorem \ref{th.SFmG} as
in the proof of Theorem \ref{thm.sine}, except that
we are now working with the squares.
%so we obtain the left-hand side of the (\ref{eq1.main_inequality}).
For the right-hand side, %of (\ref{eq1.main_inequality}),
the pinching Theorem \ref{thm.pinching} gives
$$
S\left(P_{\Z^\perp}(P_\X - P_\Y)P_{\Z^\perp} - P_\Z(P_\X - P_\Y)P_\Z \right)
\prec_w
S\left(P_\X - P_\Y\right)
$$
and we use Theorem \ref{th.SFmG} again to characterize
$S\left(P_\X - P_\Y\right)$ the same way as
in the proof of Theorem \ref{thm.sine}.
%Finally, the sine squared estimate (\ref{eq1.main_inequality})
%is equivalent to the cosine squared estimate (\ref{eq2.main_inequality})
%since evidently
%$$
%\left(|\cos^2\Theta(\X,\Z)-\cos^2\Theta(\Y,\Z))|\right)^\downarrow  =
%\left(|\sin^2\Theta(\X,\Z)-\sin^2\Theta(\Y,\Z))| \right)^\downarrow.
%$$
\end{proof}

\section{Changes in the Trial Subspace in the Rayleigh--Ritz Method}\label{section_ritz}
In this section, we explore a simple, but deep, connection
between the principal angles and the Rayleigh--Ritz method that we discuss in the Introduction.
We demonstrate that
the analysis of the influence of changes in a trial subspace in the Rayleigh--Ritz method
is a natural extension of the theory concerning principal angles
and the proximity of two subspaces developed in the previous section.

For the reader's convenience, let us repeat here
the definition of Ritz values from the Introduction:
Let $A:\H \to\H$ be a Hermitian operator and
$P_{\X}$ be an orthogonal projector to a subspace ${\X}$ of $\H$.
The eigenvalues $\Lambda(P_{\X}A|_{\X})$ are the Ritz values
of operator $A$ with respect to $\X$,
which is called the trial subspace.

Let ${\X}$ and ${\Y}$ both be subspaces of
$\H$ and $\dim{\X}=\dim{\Y}$. The goal of this
section is to analyze the sensitivity of Ritz values
with respect to the trial subspaces, specifically,
to bound the change
$\left| \Lambda\left(P_{\X}A|_{\X}\right)
-  \Lambda\left(P_{\Y}A|_{\Y}\right) \right|$
in terms of $\sin\Theta(\X,\Y)$
using weak majorization. Such an estimate
is already obtained in Theorem 10 of
our earlier paper \citet{ka03}  by
applying  Corollary \ref{cor.e_dif_majorization} to
the matrices of $P_{\X}A|_{\X}$ and $P_{\Y}A|_{\Y}$.
This approach, however, leads to an extra factor $\sqrt{2}$ on the right-hand side,
which is conjectured in \citet{ka03} to be artificial.

We remove this $\sqrt{2}$ factor in our new Theorem \ref{thm_ritz}
by using the entirely different and novel approach:
we connect the Ritz values with extension
Theorems \ref{thm.ext_to_projector} and \ref{thm.A_to_angles} on the one hand,
and with the cosine squared of principal angles in Lemma \ref{lem_projector_angles_pre}
on the other hand.
We have shown in Theorem \ref{thm.ext_to_projector}
that a Hermitian nonnegative definite contraction operator can be extended to
an orthogonal projector in a larger space. The extension
has an extra nice property: it preserves the Ritz values.
\begin{corollary}\label{cor.ritz.projtomatrix}
Under the assumptions of Theorem \ref{thm.ext_to_projector},
the Ritz values of operator $A:\H \to\H$ in the trial subspace $\X \subset H$
are the same as
the Ritz values of operator $\hat{A}:\H^2\to\H^2$ in the trial subspace
$$
\hat{\X}=\left[\begin{array}{c}
{\X} \\
0 \end{array}\right]
\subset
\hat{\H} =
\left[\begin{array}{c}
{\H} \\
0 \end{array}\right]
\subset\H^2.
$$
\end{corollary}
\begin{proof}
Let $P_{\hat{\H}}:\H^2\to\H^2$ be an orthogonal
projector on the subspace $\hat{\H}$
and $P_{\hat{\X}}:\H^2\to\H^2$ be an orthogonal
projector on the subspace $\hat{\X}$.
We use the equality sign to denote
the trivial isomorphism between $\H$ and  $\hat{\H}$,
i.e. we simply write $\H=\hat{\H}$ and $\X=\hat{\X}$.

In this notation, we first observe that
$A=P_{\hat{\H}}\hat{A}|_{\hat{\H}}$, i.e.
the operator $A$ itself can be viewed as a result of the
Rayleigh--Ritz method applied to the operator $\hat{A}$ in the trial
subspace $\hat{\H}$. Second, we use the fact that
a recursive application of the Rayleigh--Ritz method
on a system of enclosed trial subspaces is equivalent to a
direct single application of the Rayleigh--Ritz method
to the smallest trial subspace, indeed, in our notation,
$P_{\hat{\H}} P_{\hat{\X}} = P_{\hat{\X}} P_{\hat{\H}} = P_{\hat{\X}},$
since $\hat{\X} \subset \hat{\H},$ thus
$$
P_{\X}A|_{\X}=
\left.\left( P_{\hat{\X}} P_{\hat{\H}}\hat{A}|_{\hat{\H}}\right)\right|_{\hat{\X}} =
P_{\hat{\X}}\hat{A}|_{\hat{\X}}.
$$
\end{proof}

Next we note that Lemma \ref{lem_projector_angles_pre}
states that the  Rayleigh--Ritz method applied to
an orthogonal projector produces Ritz values, which are
essentially the cosines squared of the principal angles
between the range of the projector and the trial subspace.
For the reader's convenience we reformulate Lemma \ref{lem_projector_angles_pre} here:
Let the  Rayleigh--Ritz method be applied to
$A=P_{\Z}$, where $ P_{\Z}$ is an orthogonal projector
onto a subspace $\Z$,
and let ${\X}$ be the trial subspace in the Rayleigh--Ritz method;
then the set of the Ritz values is
$
\Lambda\left( P_{\X} P_{\Z}|_{\X} \right)=[ \cos^2 \Theta(\X,\Z),0,\ldots,0]
$
with $\max\{\dim\X - \dim\Z;0\}$ extra $0$s.

Now we are ready to direct our attention to the main topic of this section: the
influence of changes in a trial subspace in the Rayleigh--Ritz method
on the Ritz values.
\begin{theorem} \label{thm_ritz}
Let $A:\H\to\H$ be Hermitian
and let ${\X}$ and ${\Y}$ both be subspaces of
$\H$ and $\dim{\X}=\dim{\Y}$.
Then
\begin{equation} \label{eq.ritz}
\left| \Lambda\left(P_{\X}A|_{\X}\right)
-  \Lambda\left(P_{\Y}A|_{\Y}\right) \right|
\prec_w (\lmax-\lmin)~\sin\Theta(\X,\Y),
\end{equation}
where $\lmin$ and $\lmax$ are the smallest and largest eigenvalues of $A$,
respectively.
\end{theorem}
\begin{proof}
We prove Theorem \ref{thm_ritz} in two steps.
First we show that we can assume that
$A$ is a nonnegative definite contraction without losing
generality. Second, under these assumptions, we extend the
operator $A$ to an orthogonal projector by
Theorem \ref{thm.ext_to_projector}
and use the facts that
such an extension does not affect the Ritz values
by Corollary \ref{cor.ritz.projtomatrix}
and that
the Ritz values of an orthogonal projector can be interpreted as
the cosines squared of principal angles between subspaces
by Lemma \ref{lem_projector_angles_pre}, thus reducing
the problem to the already established result on weak majorization
of the cosine squared Theorem \ref{thm.cosine_squared}.

We observe that the statement of the theorem
is invariant with respect to a shift and a scaling, indeed,
for real $\alpha$ and $\beta$ if the operator $A$  is
replaced with $\beta(A-\alpha)$ and $\lmin$ and $\lmax$
are correspondingly updated, both sides of
(\ref{eq.ritz}) are just multiplied by $\beta$ and
(\ref{eq.ritz}) is thus invariant with respect to
$\alpha$ and $\beta$. Choosing
$\alpha=\lmin$ and $\beta=1/(\lmax-\lmin)$,
the transformed operator $(A-\lmin)/(\lmax-\lmin)$
is Hermitian with its eigenvalues enclosed
in a segment $\left[ 0,1\right]$, thus
the statement (\ref{eq.ritz}) of the theorem can be
equivalently rewritten as
\begin{equation}\label{eqn.semidef_angles}
\left|\Lambda\left( P_{\X} {A}|_{\X} \right) -
\Lambda\left( P_{\Y} {A}|_{\Y} \right)\right|
\prec_w \sin\Theta(\X,\Y),
\end{equation}
where we from now on assume that $A$ is a nonnegative definite contraction
without losing generality.

The second step of the proof is
to recast the problem into an equivalent problem
for an orthogonal projector with the same Ritz values and principal angles.
By Theorem \ref{thm.ext_to_projector} we can extend
the nonnegative definite contraction  $A$ to an orthogonal projector
$P_{\hat{\Z}}$, where $\hat{\Z}$ is a subspace of $\H^2$.
$P_{\hat{\Z}}$ has by Corollary \ref{cor.ritz.projtomatrix}
the same Ritz values with respect to trial subspaces
$$
\hat{\X}=\left[\begin{array}{c}
{\X} \\
0 \end{array}\right]
\subset
\hat{\H} =
\left[\begin{array}{c}
{\H} \\
0 \end{array}\right]
\subset\H^2
\mbox{  and   }
\hat{\Y}=\left[\begin{array}{c}
{\Y} \\
0 \end{array}\right]
\subset
\hat{\H} =
\left[\begin{array}{c}
{\H} \\
0 \end{array}\right]
\subset\H^2
$$
as $A$ has with respect to
the trial subspaces $\X$ and $\Y$.
By Lemma \ref{lem_projector_angles_pre},
these Ritz values are equal to the cosines squared of the principal angles
between $\hat{\Z}$ and the trial subspace $\hat{\X}$ or $\hat{\Y}$
possibly with the same number of $0$s being added.
Moreover, the principal angles
between $\hat{\X}$ and $\hat{\Y}$ in
$\H^2$ are clearly the same as those between
$\X$ and $\Y$ in $\H$ and
$\dim\hat{\X}=\dim\X=\dim\Y=\dim\hat{\Y}$.
Thus, (\ref{eqn.semidef_angles})
can be equivalently reformulated as
\begin{equation}\label{eqn.ritz_cosine}
|\cos^2\Theta(\hat{\X},\hat{\Z})-\cos^2\Theta(\hat{\Y},\hat{\Z}))|
\prec_w \sin\Theta(\hat{\X},\hat{\Y}).
\end{equation}
Finally, we notice that (\ref{eqn.ritz_cosine})
is already proved in
Theorem \ref{thm.cosine_squared}.
%which completes the proof of Theorem \ref{thm_ritz}.
\end{proof}

\begin{myremark} \label{r.bc}
As in \citet[Remark 7]{ka03}, the constant
$\lmax-\lmin$ in Theorem \ref{thm_ritz}
can be replaced with
$$
\max_{x \in \X+\Y,\, \|x\|=1} (x,Ax) -
\min_{x \in \X+\Y,\, \|x\|=1} (x,Ax),
$$
which for some subspaces $\X$ and $\Y$
can provide a significant improvement.
\end{myremark}

\begin{myremark}
The implications of the weak majorization inequality in Theorem \ref{thm_ritz}
may not be obvious to every reader.
To clarify, let $m=\dim{\X}=\dim{\Y}$ and
let $\alpha_1 \geq \cdots \geq \alpha_m$
be the Ritz values of $A$ with respect to ${\X}$
and $\beta_1 \geq \cdots \geq \beta_m$
be the Ritz values of $A$ with respect to ${\Y}$.
The weak majorization inequality in Theorem \ref{thm_ritz}
directly implies
$$
\sum_{i=1}^k |\alpha_i-\beta_i|^\downarrow
\leq (\lmax-\lmin) \sum_{i=1}^k \sin(\Theta_i(\X,\Y))^\downarrow
,\quad \ k=1,\ldots,m,
$$
e.g.,\ for $k=m$ we obtain
\begin{equation} \label{eq.4.8}
\sum_{i=1}^m |\alpha_i-\beta_i|
\leq (\lmax-\lmin) \sum_{i=1}^m \sin(\Theta_i(\X,\Y)),
\end{equation}
and for $k=1$ we have
\begin{equation} \label{eq.4.5}
\max_{ j=1,\ldots,m }|\alpha_j-\beta_j|
\leq (\lmax-\lmin) {\rm gap}(\X,\Y),
\end{equation}
where the gap ${\rm gap}(\X,\Y)$
between equidimensional subspaces $\X$ and $\Y$
is the sine of the largest angle between $\X$ and $\Y$.
Inequality (\ref{eq.4.5}) is proved in \citet[Th. 5]{ka03}.

For real vectors $x$ and $y$ the weak majorization $x \prec_w y$
is equivalent to the inequality
$\sum _{i=1}^n \phi(x_i) \leq \sum_{i=1}^n \phi(y_i)$
for any continuous nondecreasing convex real valued function $\phi$,
e.g.,\ \citet[Prop. 4.B.2., p. 109]{mo}. Taking, e.g., $\phi(t)=t^p$ with $p \geq 1$,
Theorem \ref{thm_ritz} also implies
$$
\Big(\sum_{i=1}^m |\alpha_i-\beta_i|^p\Big)^{\frac{1}{p}}
\leq (\lmax-\lmin) \Big(\sum_{i=1}^m \sin(\Theta_i(\X,\Y))^p\Big)^{\frac{1}{p}}
,\quad \ 1 \leq p <\infty.
$$
\end{myremark}

We finally note that the results of Theorem \ref{thm_ritz}
are not intended for the case where one of the
subspaces ${\X}$ or ${\Y}$ is invariant with respect to
operator $A.$ In such a case, it is natural to expect
a much better bound that involves the square of the
$\sin\Theta(\X,\Y)$. Majorization results of this kind
are not apparently known in the literature. Without majorization,
estimates just for the largest change in the
Ritz values are available, e.g.,\ \citet{ka03,ko05}.

\section{Application of the Majorization Results
to Graph Spectra Comparison} \label{s.g}
In this section, we show that our majorization results
can be applied to compare graph spectra.
The graph spectra comparison
can be used for graph matching and
has applications in data mining, cf. \citet{671411}.

The section is divided into three subsections.
In Subsection \ref{ss.L}, we give all
necessary definitions and basic facts
concerning Laplacian graph spectra.
In Subsection \ref{ss.gr}, we connect the Laplacian graph spectrum and Ritz values,
by introducing the graph edge Laplacian.
Finally, in Subsection \ref{ss.gm}, we prove our
main result on the Laplacian graph spectra comparison.

\subsection{Incidence Matrices and Graph Laplacians}
\label{ss.L}
Here, we give mostly well known relevant definitions,
e.g.,\ \citet{MR1324340,MR1421568,MR1275613,MR1374468,MR1637359,MR1468791},
just slightly tailored for our specific needs.

Let $V$ be a finite ordered set (of vertices),
with an individual element (vertex) denoted by $v_i \in V$.
Let $E_c$ be the finite ordered set (of all edges),
with an individual element (edge) denoted by $e_k \in E_c$
such that % $E_c \subset H^2$ and
every $e_k = [v_i,v_j]$ for all possible $i>j$. $E_c$ can be viewed as the set
of edges of a complete simple graph with vertices $V$
(without self-loops and/or multiple edges). The results
of the present paper are invariant with respect to
specific ordering of vertices and edges.

Let $w_c: E_c \to \RR$ be a function describing edge weights,
i.e. $w_c(e_k) \in \RR$. If for some edge $e_k$
the weight is positive, $w_c(e_k)>0$, we call this edge
present, %if $w_c(e_k)< 0$ we call the edge non-oriented,
if  $w_c(e_k)= 0$ we say that the edge is absent.
%Graphs with both oriented and non-oriented edges
%present are called mixed, e.g.,\ \citet{MR1955457}.
In this paper we do not allow negative edge weights.
For a given weight function $w_c$, we define
$E \subseteq E_c$ such that $e_k \in E$ if
$w_c(e_k)\neq 0$ and we define $w$ to be
the restriction of $w_c$ on all present edges $E$, i.e.
 $w$ is made of all nonzero values of $w_c$.
A pair of sets of vertices $V$
and present edges $E$ with weights $w$ is called
a graph $(V,E)$ or a weighted graph $(V,E,w)$.

The vertex--edge incidence matrix $Q_c$ of a complete graph $(V,E_c)$ is a matrix
which has a row for each vertex and a column for each edge,
with column-wise entries determined as
$q_{ik}=1,\, q_{jk}=-1$ for every edge
 $e_k = [v_i,v_j],\, i>j$ in $E_c$ and with all other entries of
$Q_c$ equal to zero.
The vertex--edge incidence matrix $Q$ of a graph $(V,E)$
is determined in the same way, but only for the
edges present in $E$. The vertex--edge incidence matrix
can be viewed as a matrix representation of a graph analog
of the divergence operator from partial differential equations
(PDE).

Extending the analogy with PDE, the matrix $L=QQ^\ast$ is called
the graph Laplacian.
%The operator form
%of the graph Laplacian $\Delta$ acting on real valued functions
%$f(v), v\in V$ is determined by
%$\Delta f (v_i) = \sum_{j: v_j\sim v_i} [f(v_i)-f(v_j)],$
%where the sum runs over the vertices joined with $v_i$
%and no boundary conditions are imposed.
In the PDE context,
this definition corresponds to the negative Laplacian
with the natural boundary conditions, cf. \citet{MR1926852}.
Let us note that in the graph theory literature
such a definitions of the graph Laplacian is usually
attributed to directed graphs, even though changing
any edge direction into the opposite does not affect
the graph Laplacian.

If we want to take into account the weights, we can
work with the matrix $Q\,\diag(w(E)) Q^\ast,$
which is an analog of an isotropic diffusion operator,
or we can introduce a more general
edge matrix $W$ and work with $Q W Q^\ast,$
which corresponds to a general anisotropic diffusion.
It is interesting to notice the equality
\begin{equation} \label{e.e}
Q_c\,\diag(w_c(E_c)) Q_c^\ast = Q\,\diag(w(E)) Q^\ast,
\end{equation}
which shows two alternative equivalent formulas for the
graph  diffusion operator.

For simplicity of presentation, we assume in the
rest of the paper that the weights $w_c$ take only the values zero
and one.
Under this assumption, we introduce matrix
$P=\diag(w_c(E_c))$   and notice that
$P$ is the matrix of an orthogonal projector
on a subspace spanned by coordinate vectors with
indices corresponding to the indices of edges present in
$E$  and that equality (\ref{e.e}) turns into
\begin{equation} \label{e.ee}
Q_c P Q_c^\ast = QQ^\ast=L.
\end{equation}
Let us note that our results can be easily extended
to a more general case of arbitrary nonnegative weights, or even to
the case of the  edge matrix $W$, assuming that
it is symmetric nonnegative definite, $W=W^\ast \geq 0.$

Fiedler's pioneering work \cite{MR0318007} on using the eigenpairs of the
graph Laplacian to determine some structural
properties of the graph
has attracted much attention in the past. Recent
advances in large-scale eigenvalue computations
using multilevel preconditioning, e.g.,
\citet{knyazev4,kn03,MR2029596},
suggest novel efficient numerical methods to
compute the Fiedler vector and
may rejuvenate this classical approach,
e.g.,\ for graph partitioning.
In this paper, we concentrate on the whole set of eigenvalues
of $L$, which is called the Laplacian graph spectrum.

It is known that the Laplacian graph spectrum does not determine the graph
uniquely, i.e. that there exist isospectral graphs,
see, e.g.,\ \citet{MR2022290} and references there.
However, intuition suggests that
a small change in a large graph should not change the
Laplacian graph spectrum very much; and attempts have been made to use
the closeness of Laplacian graph spectra to judge the
closeness of the graphs in applications;
for alternative approaches, see \citet{MR2124680}.
The goal of this section is to backup this intuition with
rigorous estimates for proximity of the Laplacian graph spectra.

\subsection{Laplacian graph spectrum and Ritz values} \label{ss.gr}
In the previous section, we obtain in Theorem \ref{thm_ritz}
a weak majorization bound for changes in Ritz values depending on a change
in the trial subspace, which we would like to apply to
analyze the graph spectrum. In this subsection,
we present an approach that
allows us to interpret the Laplacian graph spectrum
as a set of Ritz values obtained by the
Rayleigh--Ritz method applied to the complete graph.

A graph $(V,E)$ can evidently be obtained from the
complete graph $(V,E_c)$ by removing edges, moreover,
as we already discussed, we can construct the
 $(V,E)$ graph Laplacian by either of the terms
 in equality (\ref{e.ee}).
The problem is that such a construction cannot be
recast as an application of the Rayleigh--Ritz method,
since the multiplication by the projector $P$
takes place inside of the product in (\ref{e.ee}), not outside, as
required by the  Rayleigh--Ritz method.

To resolve this difficulty, we use
the matrix $K=Q^\ast Q$ that is sometimes called
the matrix of the graph {\em edge} Laplacian, instead of
the matrix of the graph {\em vertex} Laplacian $L=QQ^\ast,$
as both matrices $K$ and $L$ share the same nonzero eigenvalues.
The advantage of the edge Laplacian $K$
is that it can be obtained from the
edge Laplacian of the complete graph
$Q_c^\ast Q_c$
simply by removing the rows and columns
that correspond to missing edges. Mathematically,
this procedure can be viewed as an instance of the classical Rayleigh--Ritz method:
\begin{lemma} \label{thm.K}
Let  us remind the reader that
the weights $w_c$ take only the values zero and one
and that $P=\diag(w_c(E_c))$ is a matrix of an orthogonal projector
on a subspace spanned by coordinate vectors with
indices corresponding to the indices of edges present in $E$. Then
$
Q^\ast Q= \left.\left( P Q_c^\ast Q_c \right)\right|_{{\rm Range}(P)},
$
in other words, the matrix $Q^\ast Q$ is the result
of the Rayleigh--Ritz method applied to the matrix
$Q_c^\ast Q_c$ on the trial subspace ${\rm Range}(P)$.
\end{lemma}
The application of the Rayleigh--Ritz method in this case
is reduced to simply crossing out
rows and columns of the matrix $Q_c^\ast Q_c$
corresponding to absent edges,
since $P$ projects onto a span of coordinate vectors
with the indices of the present edges.

Lemma \ref{thm.K} is a standard tool in
the spectral graph theory, e.g.,\ \citet{MR1344588},
to prove the eigenvalues interlacing; however,
the procedure is not apparently recognized in the spectral graph community
as an instance of the classical Rayleigh--Ritz method.
Lemma \ref{thm.K} provides us with the missing link
in order to apply our Theorem \ref{thm_ritz} to
Laplacian graph spectra comparison.

\subsection{Majorization of Ritz Values for Laplacian Graph Spectra Comparison}
\label{ss.gm}
Using the tools that we have presented in the previous subsections,
we now can apply the particular case, (\ref{eq.4.8}), of our weak majorization
result of Section \ref{section_ritz}
to analyze the change in the graph spectrum when several edges are
added to or removed from the graph.

\begin{theorem} \label{thm.g}
Let $(V,E^1)$ and $(V,E^2)$ be two graphs
with the same set of $n$ vertices $V$, with the same
number of edges $E^1$ and $E^2,$ and with the number of
differing edges in $E^1$ and $E^2$ equal to $l$.
Then
\begin{equation} \label{e.g}
\sum_k |\lambda^1_k - \lambda^2_k| \leq n l,
\end{equation}
where $\lambda^1_k$ and $\lambda^2_k$ are all
elements of the Laplacian spectra of the graphs
$(V,E^1)$ and $(V,E^2)$ in nonincreasing order.
\end{theorem}
\begin{proof}
The spectra of the graph vertex and edge Laplacians $QQ^\ast$ and
$Q^\ast Q$ are the same apart from zero, which does not affect the
statement of the theorem, so we redefine $\lambda^1_k$ and
$\lambda^2_k$ as elements of the spectra, counting the multiplicities,
of the edge Laplacians of
the graphs $(V,E^1)$ and $(V,E^2)$. Then, by Theorem \ref{thm.K},
$\lambda^1_k$ and $\lambda^2_k$ are the Ritz values of the edge
Laplacian matrix $A=Q_c^\ast Q_c$ of the complete graph,
corresponding to the trial subspaces $\X={\rm Range}(P_1)$ and
$\Y={\rm Range}(P_2)$ spanned by coordinate vectors with indices
of the edges present in  $E^1$ and $E^2$, respectively.

Let us apply Theorem \ref{thm_ritz}, taking the sum
over all available nonzero values in the weak majorization statement
as in  (\ref{eq.4.8}).
This already gives us the left-hand side of (\ref{e.g}).
To obtain the right-hand side of (\ref{e.g}) from Theorem \ref{thm_ritz},
we now show in our case that, first, $\lmax-\lmin = n$ and, second,
the sum of sines of all angles between the trial subspaces
$\X$ and $\Y$ is equal to $l$.

The first claim follows from the fact, which is easy to check
by direct calculation,
that the spectrum of the vertex (and thus the edge) Laplacian
of the complete graph with $n$ vertices
consists of only two eigenvalues $\lmax= n$ and $\lmin=0$.
Let us make a side note that we can interpret
the Laplacian of the complete graph as a scaled projector, i.e.
in this case we could have applied Theorem \ref{thm.cosine_squared}
directly, rather than Theorem \ref{thm_ritz}, which would still result in
(\ref{e.g}).

The second claim, on the sum of sines of all angles, follows from the
definition of $\X$ and $\Y$ and the assumption
that the number of
differing edges in $E^1$ and $E^2$ is equal to $l$.
Indeed, $\X$ and $\Y$ are spanned by coordinate vectors
with indices of the edges present in  $E^1$ and $E^2$.
The edges that are present both in  $E^1$ and $E^2$
contribute zero angles into $\Theta(\X,\Y)$,
while the $l$ edges that are different in  $E^1$ and $E^2$
contribute $l$ right angles into $\Theta(\X,\Y)$,
so that the sum of all terms in $\sin\Theta(\X,\Y)$
is equal to $l$.
\end{proof}

Remark \ref{r.bc} is also applicable for
Theorem \ref{thm.g} --- while the $\min$ term is
always zero, since all graph Laplacians
are degenerate, the $\max$ term can be made smaller by
replacing $n$ with the largest eigenvalue
of the Laplacian of the graph $(V,E^1\cup E^2)$.

It is clear from the proof that we do not use the full force
of our weak majorization results in Theorem \ref{thm.g},
because it concerns angles which are zero or $\pi/2.$
Nevertheless, the results of Theorem \ref{thm.g}
appear to be novel in graph theory.
We note that these results can be easily extended on {\it k}-partite graphs,
and possibly to mixed graphs.

Let us finally mention an alternative approach to compare
Laplacian graph spectra, which we do not cover in the present paper,
by applying Corollary \ref{cor.e_dif_majorization} directly to
graph Laplacians and estimating the right-hand side using the
fact that the changes in $l$ edges represents a low--rank perturbation
of the graph Laplacian, cf. \citet{MR1708588}.

\section*{Conclusions}
We use majorization to investigate the sensitivity of angles
between subspaces and  Ritz values with respect to subspaces,
and to analyze changes in graph Laplacian spectra where edges are added and removed.
We discover that these seemingly different areas are all surprisingly related.
We establish in a unified way new results on weak majorization
of the changes in the sine/cosine (squared) and in the Ritz values.
The main strength of the paper in our opinion is, however,
not so much in the results themselves but rather
in a novel and elegant proof technique that is based on
a classical but rarely used idea of
extending Hermitian operators to orthogonal projectors in a larger space.
We believe that this technique is very powerful and should be
known to a wider audience.

\section*{Acknowledgments}
We are indebted to Li Qiu and Yanxia Zhang, who have suggested to us the idea
of the proof of the sine estimate in Theorem \ref{thm.sine}
that has allowed us also to simplify
our original proof of Theorem \ref{thm.cosine_squared}.
The authors thank our Ph.D. students Ilya Lashuk for contributing to Section \ref{section_extend_proj}
and Abram Jujunashvili for many helpful comments.
Finally, we thank anonymous referees for their numerous
useful and illuminating suggestions that have dramatically improved the paper.

\bibliographystyle{plainnat}
%\bibliography{argentati}

\end{document}